\documentclass[11pt]{amsart}
\usepackage{amsmath}
\usepackage{amssymb}
\usepackage{amsfonts}
\usepackage{mathrsfs}
\usepackage[a4paper]{geometry}
\theoremstyle{plain}
\newtheorem{teo}{Theorem}

\newtheorem{lem}{Lemma}

\theoremstyle{definition}

\theoremstyle{remark}

\newtheorem{assum}{Assumption}

\newcommand{\re}{\mathbb{R}}

\newcommand{\g}{\Gamma}
\newcommand{\ja}{P_n^{\alpha,\beta}}
\def\d{\displaystyle}
\def\pn{\par\noindent}

\def\videbox{\mathbin{\vbox{\hrule\hbox{\vrule height1ex \kern.5em\vrule
height1ex}\hrule}}}
\def\1{1\!{\rm l}}

\title{Large deviations for statistics of Jacobi process}
\date{\today \\ {\it Key Words}:  Jacobi Process, Subordinated Jacobi process, Large deviations, Maximum likelihood.}
\begin{document}
\maketitle
\centerline{N. DEMNI
\footnote{Laboratoire de Probabilit\'es et Mod\`eles Al\'eatoires, Universit\'e de Paris VI, 4 Place Jussieu, Case 188, F-75252 Paris Cedex 05.}
AND M. ZANI \footnote{Laboratoire d'analyse et de math\'ematiques appliqu\'ees, Universit\'e de Paris XII, Val de Marne.}}
\vspace{1cm}
\begin{abstract} This paper is aimed to derive large deviations for statistics of Jacobi process already conjectured by M. Zani in her Thesis. To proceed, we write in a more simple way the Jacobi semi-group density. Being given by a bilinear sum involving Jacobi polynomials, it differs from Hermite and Laguerre cases by the quadratic form of its eigenvalues. Our attempt relies on subordinating  the process using a suitable random time-change. This will give an analogue of Mehler formula whence we can recover the desired expression by inverting some Laplace transforms. Once we did, an adaptation of Zani's result in the non-steepness case will provide the required large deviations principle. 
\end{abstract}

\section{Introduction} The Jacobi process is a Markov process on $[-1,1]$ given by the following infinitesimal generator:
\begin{equation*}\mathscr{L} = (1-x^2)\frac{\partial^2}{\partial^2 x} + (px+q)\frac{\partial}{\partial x}, \qquad x \in [-1, 1]\end{equation*} for some real $p, q$, defined up to the first time when it hits the boundary.
In fact, it belongs to the class of diffusions associated to some families of orthogonal polynomials, i.e. the
infinitesmal generator admits an orthogonal polynomials basis as  eigenfunctions (\cite{Ba}) such as Hermite, Laguerre
and Jacobi polynomials . More precisely, if  $\ja$ denotes  the Jacobi polynomial with parameters $\alpha, \beta > -1$
defined by :\begin{equation*}
\ja(x) = \frac{(\alpha+1)_n}{n!}{}_2F_1\left(-n,n+\alpha+\beta+1,\alpha+1; \frac{1-x}{2}\right), \qquad x \in [-1,1],
\end{equation*} then we can see that :
\begin{equation*}\mathscr{L}\ja = -n(n+\alpha + \beta +1)\ja
\end{equation*} 
for $p= - (\beta + \alpha +2)$ and $q=\beta - \alpha$. In 1964, Wong resolved the forward Kolmogorov or Fokker-Planck equation  (see \cite{Wong}, \cite{Rev})
\begin{equation*}
\partial_y^2[B(y)p] - \partial_y[A(y)p] = \partial_t p, \quad p = p_t(x,y),
\end{equation*}
where $B,\,A$ are polynomials of degree $2,\,1$ respectively, and gave the principal solution $(p_0(x,y) = \delta_x(y))$ using the classical Sturm-Liouville theory. This gives rise to 
a class of stationnary Markov processes satisfying :
\begin{equation*}
\lim_{t \rightarrow \infty}p_t(x,y) = \int_{x_1}^{x_2} W(x)p_t(x,y) dx = W(y)
\end{equation*}
where $W$ is the density function solution of the corresponding Pearson equation (\cite{Wong}). In our case, $p_t$ has the discrete spectral decomposition :
\begin{equation}p_t(x, y) = \left(\sum_{n\geq 0}(R_n)^{-1}e^{-\lambda_nt}\ja(x)\ja(y)\right)W(y), \qquad x, y \in [-1, 1]
\end{equation}where
\begin{equation*}
\lambda_n = n(n + \alpha + \beta +1), \qquad W(y) = \frac{(1-y)^{\alpha}(1+y)^{\beta}}{2^{\alpha+\beta+1}B(\alpha+1,\beta+1)}
\end{equation*}
with $B$ denoting the Beta function and (\cite{Askey}, p. 99) : 
\begin{equation*}
R_n\footnote{$(\ja(x))_{n \geq 0}$ are normalized such that they form an orthogonal basis with respect to the probability measure $W(y)dy$ which is not the same used in \cite{Askey}.} 
= ||\ja||^2_{L^2([-1, 1] , W(y)dy)} = \frac{\g(\alpha+\beta+2)}{2n + \alpha + \beta +1}\frac{(\alpha+1)_n(\beta+1)_n}{\g(\alpha+\beta+n+1)n!}
\end{equation*}
 
Few years later, Gasper (\cite{Gas}) showed that this bilinear sum is the transition kernel of a diffusion and that is a solution of the heat equation governed by a Jacobi operator, generalizing a previous result of Bochner for ultraspherical polynomials (\cite{Boc}). It is worthnoting that $\lambda_n$ has a quadratic form while in the Hermite (Brownian) and Laguerre (squared Bessel) cases $\lambda_n = n$.
Hence, we will try to subordinate the Jacobi process by the mean of a random time-change in order to get a Mehler type formula. What is quite interesting is that subordinated Jacobi process semi-group, say $q_t(x, y)$, is the Laplace transform of $p_{2/t}(x, y)$. Thus, we recover a suitable expression for $p_t(x, y)$ by inverting some Laplace transforms already computed by Biane, Pitman and Yor (see \cite{Biane}, \cite{Pit}). This will allows us to derive a LDP for the maximum likelihood estimate (MLE) for $p$ in the ultraspherical case, i. e. $q=0 \,  (\beta = \alpha)$, a fact conjectured by Zani in her thesis. Then, using a skew product representation of the Jacobi process involving squared Bessel processes, we construct a
family $\{\hat{\nu}_t\}_t$ of estimators for the index $\nu$ of the squared Bessel process based on a Jacobi trajectory observed till time $t$. This satisfies a LDP with the same rate function derived for the MLE based on a squared Bessel trajectory. 

\subsection{Inverse Gaussian subordinator}
By an {\it inverse Gaussian subordinator}, we mean the process of the first hitting time of a Brownian motion with drift $B_t^{\mu}:= B_t + \mu t$, $\mu > 0$, namely,\begin{equation*}
T_t^{\mu,\delta} = \inf\{s>0 ; \quad B_s^{\mu} = \delta t\}, \qquad \delta > 0.\end{equation*}
Using martingale methods, we can show that for each $t > 0$, $u \geq 0$,
\begin{equation*}
\mathbb{E}(e^{-uT_t^{\mu,\delta}}) = e^{-t\delta(\sqrt{2u + \mu^2} - \mu)}
\end{equation*}
whence we recover the density below :
\begin{equation*}
\nu_t(s) = \frac{\delta t}{\sqrt{2\pi}}e^{\delta t\mu } s^{-3/2}\exp\left(-\frac{1}{2}(\frac{t^2\delta^2}{s} + \mu^2s)\right){\bf 1}_{\{s>0\}}
\end{equation*}

\subsection{The subordinated Jacobi Process}
Let us consider a Jacobi process $(X_t)_{t \geq 0}$. Then, using $(1)$, the semi-group of the subordinated Jacobi process $(X_{T_t^{\mu,\delta}})_{t \geq 0}$ is given by:
\begin{align*}
q_t(x,y) &= \int_0^{\infty} p_s(x,y) \nu_t(s)ds \\&= W(y)\sum_{n\geq 0}(R_n)^{-1}\left(\int_0^{\infty}e^{-\lambda_ns}\nu_t(s)ds\right)\ja(x)\ja(y)
\\&= W(y)\sum_{n\geq 0}(R_n)^{-1}\mathbb{E}(e^{-\lambda_nT_t^{\mu,\delta}})\ja(x)\ja(y)
\end{align*}
Writing $\lambda_n = (n + \gamma)^2 - {\gamma}^2$ where $\gamma = \displaystyle \frac{\alpha + \beta +1}{2}$, and substituting $\delta = 1/\sqrt 2$, $\mu = {\sqrt 2}\gamma$
for $\alpha + \beta > -1$ in the expression of $\nu_t$, one gets :
\begin{equation*}
\mathbb{E}(e^{-\lambda_nT_t^{\mu,\delta}}) = e^{-nt}
\end{equation*}
so that
\begin{equation*}
q_t(x, y) =  W(y)\sum_{n\geq 0}(R_n)^{-1}e^{-nt}\ja(x)\ja(y)
\end{equation*}
The last sum has been already computed (\cite{Askey}, p. 385) :
\begin{align}
\sum_{n=0}^{\infty}(R_n)^{-1}\ja(x)\ja(y)r^n &= \frac{1-r}{(1+r)^a}
\sum_{m,n \geq 0}\frac{\left(\frac{a}{2}\right)_{m+n}\left(\frac{a+1}{2}\right)_{m+n}}{(\alpha+1)_m(\beta+1)_n}\frac{u^mv^n}{m!n!}  \nonumber
\\& = \frac{1-r}{(1+r)^a} \, F_4(\frac{a}{2}, \frac{a+1}{2}, \alpha +1, \beta +1; u, v) \label{A}
\end{align}  where $|r| < 1$, $a = \alpha + \beta +2$, $F_4$ is the Appell function and
\begin{equation*}
u = \frac{(1-x)(1-y)r}{(1+r)^2} \qquad v =  \frac{(1+x)(1+y)r}{(1+r)^2}.
\end{equation*} 
Then, we use the integral representation of $F_4$ (see \cite{Exton}, p 51) to get:
\begin{align*}
q_t(x, y) &= \frac{W(y)}{\g(a)}\frac{1-r}{(1+r)^a}  \int_0^{\infty}s^{a-1}e^{-s}{}_0F_1(\alpha+1; \frac{u}{4}s^2){}_0F_1(\beta+1; \frac{v}{4}s^2)ds
\\& \overset{(1)}{=} \frac{W(y)}{\g(a)}\frac{1-r}{(1+r)^a}\int_0^{\infty}s^{a-1}e^{-s}\sum_{n\geq0}\frac{\ja(z)}{(\alpha+1)_n(\beta+1)_n} A^ns^{2n} \,ds
\\& \overset{(2)}{=} \frac{W(y)}{\g(a)}\frac{1-r}{(1+r)^a}\sum_{n\geq0} \frac{\g(2n+a)}{(\alpha+1)_n(\beta+1)_n} \ja(z)A^n
\end{align*}
where in $(1)$, we used (see \cite{Mag}, p 214)
\begin{equation*}
{}_0F_1(c; w(1-r)/2){}_0F_1(d; w(1+r)/2) = \sum_{n\geq0}\frac{\ja(r)}{(c)_n(d)_n} w^n, \, \alpha = c-1, \, \beta = d-1,
\end{equation*}
 in $(2)$, we used Fubini Theorem, $z = \displaystyle \frac{x+y}{1+xy}$ and $A= \displaystyle
\frac{(1+xy)r}{2(1+r)^2}$. Letting $r=e^{-t}$, then
\begin{align*}
q_t&(x, y) = \frac{W(y)e^{\frac{a-1}{2}t}}{2^{a-1}}\frac{\sinh(t/2)}{(\cosh(t/2))^a}
\sum_{n\geq0} \frac{(a)_{2n}}{(\alpha+1)_n(\beta+1)_n} \ja(z)\left[\frac{(1+xy)}{8\cosh^2(t/2)}\right]^n
\\&=  
\frac{W(y)\tanh(t/2)e^{\frac{a-1}{2}t}}{2^{a-1}}
\sum_{n\geq0} \frac{(a)_{2n}}{(\alpha+1)_n(\beta+1)_n} \ja(z)\left[\frac{(1+xy)}{8}\right]^n \left(\frac{1}{\cosh(t/2)}\right)^{2n+a-1}.
\end{align*}
Besides,
\begin{align*}
q_t(x, y) = \frac{t\, e^{\gamma t}}{2\sqrt{\pi}}\int_0^{\infty}p_s(x, y)\, s^{-3/2}e^{-{\gamma}^2s}\,e^{-\frac{t^2}{4s}}ds
 = \frac{t\, e^{\gamma t}}{2\sqrt{2\pi}}\int_0^{\infty}p_{2/r}(x, y)\, r^{-1/2}e^{-2{\gamma}^2/r}\,e^{-\frac{t^2}{8}r}dr
 \end{align*}
Thus, noting that $\gamma = (a-1)/2$, we get :
\begin{align*}\int_0^{\infty} & p_{2/r}(x, y)\, r^{-1/2}e^{-2{\gamma}^2/r}\,e^{-\frac{t^2}{8}r}dr =  \frac{\sqrt{2\pi}W(y)}{2^{a-1}}\frac{\tanh(t/2)}{t/2}
\\&  \sum_{n\geq0} \frac{(a)_{2n}}{(\alpha+1)_n(\beta+1)_n} \ja(z)\left[\frac{(1+xy)}{8}\right]^n \left(\frac{1}{\cosh(t/2)}\right)^{2n+a-1}.
\end{align*}

\subsection{The Jacobi semi-group}
The following results are due to Biane, Pitman and Yor (see \cite{Biane}, \cite{Pit}) :
\begin{eqnarray}
\label{e3} \int_0^{\infty}e^{-\frac{t^2}{8}s}f_{C_h}(s)\,ds &=& \left(\frac{1}{\cosh(t/2)}\right)^h, \qquad h>0 \\
\label{e4} \int_0^{\infty}e^{-\frac{t^2}{8}s}f_{T_h}(s)\,ds &=& \left(\frac{\tanh(t/2)}{(t/2)}\right)^h , \qquad h>0
\end{eqnarray}
where $(C_h)$ and $(T_h)$ are two families of L\'evy processes with respective density functions $f_{C_h}$ and $f_{T_h}$ for fixed $h > 0$.
The densities of $C_h$ and $T_1$ are given by (\cite{Biane}):
\begin{eqnarray*}
f_{C_h}(s) & = & \frac{2^h}{\g(h)}\sum_{p \geq 0}(-1)^p\frac{\g(p+h)}{p!}f_{\tau(2p+h)}(s)\\
f_{T_1}(s) & = & \sum_{k \geq 0}e^{-\frac{{\pi}^2}{2}(k+\frac{1}{2})^2s}{\bf 1}_{\{s >0\}}
\end{eqnarray*} where $\tau(c) = \inf\{r>0 ; B_r = c\}$ is the L\'evy subordinator ( i. e, the first hitting time of a standard Brownian motion $B$) with corresponding density :
\begin{equation*}
f_{\tau(2p+h)}(s) = \frac{(2p+h)}{\sqrt{2\pi s^3}}\exp\left(-\frac{(2p+h)^2}{2s}\right) {\bf 1}_{\{s > 0\}}.
\end{equation*}
Thus :
\begin{align*}
p_{2/r}(x, y) =  \frac{\sqrt{2\pi r}W(y)e^{2{\gamma}^2/r}}{2^{a-1}}
\sum_{n\geq0} \frac{(a)_{2n}}{(\alpha+1)_n(\beta+1)_n} \ja(z) \left[\frac{(1+xy)}{8}\right]^n  \times \left(f_{T_1} \star f_{C_{2n+a-1}}\right)(r)
\end{align*}
or equivalently (where $B$ stands for the Beta function) :
\begin{align*}
p_{t}(x, y) =  \frac{\sqrt{\pi}W(y)}{2^{\alpha+\beta}}\frac{e^{\gamma^2t}}{\sqrt{t}} \sum_{n\geq0} \frac{(a)_{2n}}{(\alpha+1)_n(\beta+1)_n}
\ja(z)\left[\frac{(1+xy)}{8}\right]^n  \left(f_{T_1} \star f_{C_{2n+\alpha+\beta+1}}\right)(\frac{2}{t}).
 \end{align*}

\subsection{The ultraspherical case}
This case corresponds to $\alpha = \beta > \displaystyle -\frac{1}{2}$ and we will proceed slight differently. Indeed, $a = 2\alpha +2$ and 
\begin{align*}
(\ref{A}) &= \frac{1-r}{(1+r)^{2\alpha+2}} \, F_4(\alpha +1, \alpha + 3/2, \alpha +1, \alpha +1; u, v)
\\& =\frac{1-r}{(1+r)^{2\alpha+2}} \, \frac{1}{(1-u-v)^{\alpha+3/2}}{}_2F_1(\frac{2\alpha+3}{4}, \frac{2\alpha+5}{4}, \alpha + 1; \frac{4uv}{(1-u-v)^2})
\end{align*}
where the last equality follows from (see \cite{Bry}) 
\begin{equation*}
F_4(b, c, b, b; u, v) = (1-u-v)^{-c}{}_2F_1(c/2, (c+1)/2, b; \frac{4uv}{(1-u-v)^2}).
\end{equation*} 
 Hence,
\begin{align*}
&q_t(x, y) = \frac{W(y)e^{\frac{2\alpha+1}{2}t}}{2^{\alpha+1/2}} \frac{\sinh(t)}{(\cosh t -xy)^{\alpha+3/2}}
{}_2F_1(\frac{2\alpha+3}{4}, \frac{2\alpha+5}{4}, \alpha + 1; \frac{(1-x^2)(1-y^2)}{(\cosh t -xy)^2})
\\& = \frac{W(y)e^{\frac{2\alpha+1}{2}t}}{2^{\alpha+ 1/2}} \sinh(t)
\sum_{n\geq 0}\frac{[(2\alpha+3)/4]_n[(2\alpha+5)/4]_n}{(\alpha + 1)_n}\frac{[(1-x^2)(1-y^2)]^n}{(\cosh t - xy)^{2n + \alpha +3/2}}.
\end{align*}
Besides, for $h>0$, we may write :
\begin{align*}
\left(\frac{1}{\cosh t - xy}\right)^h = \sum_{k \geq 0}\frac{(h)_k}{k!}\, \frac{(xy)^k}{(\cosh t)^{k+h}}
\end{align*}since $\displaystyle \left|\frac{xy}{\cosh t}\right| < 1\quad \forall x, y \in ]-1, 1[ ,\, \forall t \geq 0$ and where we used:
\begin{equation*}
\frac{1}{(1-r)^h} = \sum_{k \geq 0}\frac{(h)_k}{k!}\, r^k \qquad h > 0,\, |r| < 1.
\end{equation*}
Consequently, using Gauss duplication formula,
\begin{align*}
q_t(x, y) 
 = K_{\alpha}W(y)e^{\frac{2\alpha+1}{2}t}\tanh(t)
 \sum_{n, k \geq 0}\frac{\g(\nu(n,k,\alpha) + 1)(xy)^k}{k!\, n! \, \g(\alpha + n +1)} \left[\frac{(1-x^2)(1-y^2)}{4}\right]^n \left(\frac{1}{\cosh t}\right)^{\nu(n,k,\alpha)}
\end{align*} 
where $\nu(n,k,\alpha) = 2n + k + \alpha + 1/2$ and $K_{\alpha} = \g(\alpha+1)/[2^{\alpha+1/2}\g(\alpha+3/2)]$. 
Thus, since $\gamma = \alpha + 1/2$ when $\alpha = \beta$, one has :
\begin{align*}
\int_0^{\infty}p_s(x, y)\, s^{-3/2}e^{-{\gamma}^2s}\,e^{-\frac{t^2}{4s}}ds &= \frac{\sqrt{2\pi}\g(\alpha+1)}{2^{\alpha}\g(\alpha+3/2)} \frac{\tanh(t)}{t} W(y)
 \\& \sum_{n, k \geq 0}\frac{\g(\nu(n,k,\alpha) + 1)(xy)^k}{k!\, n! \, \g(\alpha + n +1)} \left[\frac{(1-x^2)(1-y^2)}{4}\right]^n \left(\frac{1}{\cosh t}\right)^{\nu(n,k,\alpha)}
\end{align*}
or equivalently:
\begin{align*}
\int_0^{\infty}p_{1/2s}(x, y)e^{-\frac{{\gamma^2}}{2s}} e^{-\frac{t^2}{2}s}&\frac{ds}{\sqrt s} = \frac{\sqrt{\pi}\g(\alpha+1)}{2^{\alpha}\g(\alpha+3/2)} \frac{\tanh(t)}{t} W(y)
\\& \sum_{n, k \geq 0}\frac{\g(\nu(n,k,\alpha) + 1)(xy)^k}{k!\, n! \, \g(\alpha + n +1)} \left[\frac{(1-x^2)(1-y^2)}{4}\right]^n \left(\frac{1}{\cosh t}\right)^{\nu(n,k,\alpha)}
\end{align*}
Using $(\ref{e3})$, $(\ref{e4})$, $f_{C_h}$ et $f_{T_1}$   (we take $t^2/2$ instead of $t^2/8$), the density is written :
\begin{align*}
p_{1/2s}(x, y) &= \frac{\sqrt{\pi s}\g(\alpha+1)}{2^{\alpha}\g(\alpha+3/2)}  W(y)e^{\frac{{\gamma^2}}{2s}}
\\& \sum_{n, k \geq 0}\frac{\g(\nu(n,k,\alpha)+ 1)}{k!n! \g(\alpha + n +1)} \left(\frac{xy}{2}\right)^k\left[\frac{(1-x^2)(1-y^2)}{4}\right]^n  f_{T_1} \star f_{C_{\nu(n,k,\alpha)}}(s)
\end{align*}
Finally 
\begin{align*}
p_{t}(x, y) &=  \sqrt{\pi}K_{\alpha} \frac{e^{\gamma^2t}}{\sqrt{t}}W(y)\sum_{n, k \geq 0}\frac{\g( \nu(n,k,\alpha) + 1)(xy)^k}{k!n! \g(\alpha + n +1)} 
 \left[\frac{(1-x^2)(1-y^2)}{4}\right]^nf_{T_1} \star f_{C_{\nu(n,k,\alpha)}}(\frac{1}{2t})
\end{align*}

\section{Application to statistics for diffusions processes}

\subsection{Some properties of the Jacobi process}
In the probability scope, we are used to define the Jacobi process on $[-1,1]$ as the unique strong solution of the SDE :
\begin{equation*}
dY_t = \sqrt{1-Y_t^2}dW_t + (bY_t + c) dt.
\end{equation*}
It is straightforward that $(Y_t)_{t \geq 0} \overset{\mathscr{L}}{=} (X_{t/2})_{t \geq 0}$ where $X$ is the Jacobi process already defined in section 1 with $p = 2b, \, q = 2c$. In order to derive some facts, let us make the variable change $y \mapsto (y+1)/2$, this gives up to a time change $(t \rightarrow 4t)$ :
\begin{align*}
dJ_t & = 2\sqrt{J_t(1-J_t)}dW_t + [2(c-b) + 4bJ_t]\, dt
       \\& = 2\sqrt{J_t(1-J_t)}dW_t + [d - (d + d')J_t]\, dt
\end{align*}
where $d = 2(c-b) = q-p = 2(\beta + 1)$ and $d' = -2(c+b) = -(p+q) = 2(\alpha+ 1)$, which is the Jacobi process of parameters $(d,d')$ already considered in \cite{War}.  Moreover, authors provide the following skew-product : let $Z_1, Z_2$ be two independent Bessel processes of dimensions $d, \,d'$ and starting from $z,\,z'$ respectively. Then :
\begin{equation*}
\left(\frac{Z _1^2(t)}{Z_1^2(t) + Z_2^2(t)}\right)_{t \geq 0} \overset{\mathscr{L}}{=} \left(J_{A_t}\right)_{t \geq 0}, \quad A_t := \int_0^t \frac{ds}{Z_1^2(s) + Z_2^2(s)}, \quad J_0 = \frac{z}{z+z'}.
\end{equation*}
Using well known properties of squared Bessel processes (see \cite{Rev}), one deduce that if $d \geq 2 \, (\beta \geq 0)$ and $z > 0$, then $J_t > 0$ a. s. for all $t > 0$.
Since $1-J$ is still a Jacobi process of parameters $(d',d)$, then, for $d' \geq 2, \, (\alpha \geq 0)$ and $z' > 0$, $J_t < 1$ a. s. for all  $t > 0$.
These results fit in the one dimensional case those established in \cite{Yan} for the matrix Jacobi process (Theorem 3. 3. 2, p. 36).
Besides, since $0$ is a reflecting boundary for $Z_1, \, Z_2$ when $0 < d, \, d' < 2\, (-1< \alpha,\beta < 0)$, then both $0$ and $1$ are reflecting boundaries for $J$.

\subsection{LDP in the ultraspherical case}
Let us consider the following SDE corresponding to the ultraspherical Jacobi process:

\begin{equation}\label{eds}
\begin{cases}
dY_t=\sqrt{1-Y_t^2}dW_t+bY_tdt \\
Y_0  = y_0 \in ]-1,1[ \end{cases}
\end{equation}
Let $Q_{y_0}^b$ be the law of $(Y_t,\,t\geq 0)$ on the canonical filtered probability space $(\Omega,({\mathcal F}_t),{\mathcal F})$ where $\Omega$ is the space
of $]-1,1[$--valued functions. The parameter $b$ is such that $b \leq -1$ (or $\alpha \geq 0$), so that $-1 < Y_t < -1$ for all $t > 0$.
The maximum likelihood estimate of $b$ based on the observation of a single trajectory $(Y_s,\,0\leq s\leq t)$ under $Q_0^b$ is given by
\begin{equation}\label{mle}
\hat b_t=\d\frac{\int_0^t\frac{Y_s}{1-Y_s^2}\,dY_s}{\int_0^t\frac{Y_s^2}{1-Y_s^2}\,ds}.
\end{equation}
The main result of this section is the following theorem.
%%%%%%%%%%%%%%%%%%%%%%%%%%%%%%%%%%%%%%%%%%%%%%%%%%%%%%%%%%%%%%%%%%%%%%%
\begin{teo}
\label{pgd1} When $b \leq -1$, the family $\{\hat b_{t}\}_t$ satisfies a
LDP with speed $t$ and good rate function

\begin{equation}
J_b(x)=\begin{cases} \d-\frac{1}{4}\frac{(x-b)^2}{x+1}  \quad  {\textrm if} \quad  x\leq x_0 \\
x+2+\sqrt{(b-x)^2+4(x+1)} \quad {\textrm if} \quad  x > x_0 > b
\end{cases}
\end{equation}

where $x_0$ is the unique solution $x<-1$ of the equation
$$(b-x)^2=4x(x+1)\,.$$
\end{teo}
%%%%%%%%%%%%%%%%%%%%%%%%%%%%%%%%%%%%%%%%%%%%%%%%%%%%%%%%%%%%%%%%%%%%%%%%
\medskip\pn

\medskip\pn
{\it Proof of Theorem \ref{pgd1}:}
 \pn We follow the scheme of
Theorem 3.1 in \cite{Zani}. Let us  denote by :
$$S_{t,x}=\int_0^t\frac{Y_s}{1-Y_s^2}\,dY_s-x\int_0^t\frac{Y_s^2}{1-Y_s^2}\,ds$$
so that for $x>b$ (resp. $x<b$), $P(\hat b_t\geq x)=P(S_{t,x}\geq
0)$ (resp. $P(\hat b_t\leq x)=P(S_{t,x}\leq 0)$). Therefore, to
derive a large deviation principle on $\{\hat b_t\}$, we seek a LDP
result for $S_{t,x}/t$ at $0$. Let us compute the normalized cumulant
generating function $\Lambda_{t,x}$ of $S_{t,x}$:
\begin{equation}
\Lambda_{t,x}(\phi)=\frac{1}{t}\log E(e^{\phi S_{t,x}})
\end{equation}
From Girsanov formula, the generalized densities are given by
$$\frac{dQ_a^{b}}{dQ_a^{{b_0}}}=\exp\left\{({b}-{b_0})\int_0^t\frac{Y_s}{1-Y_s^2}\,dY_s - \frac{1}{2}({b}^2-{b_0}^2)\int_0^t\frac{Y_s^2}{1-Y_s^2}\,ds \right\}
$$
From It\^o formula, $$F(Y_t)=-\frac{1}{2}\log(1-Y_t^2)=\int_0^t\frac{Y_s}{1-Y_s^2}\,dY_s  + \frac{1}{2}\int_0^t\frac{1+ Y_s^2}{1-Y_s^2}\,ds.$$
Let us denote by
$${\mathcal D}_1=\{\phi\,:(b+1)^2+2\phi(x+1)\geq 0\}\,.$$
Set $b(\phi,x)=-1-\sqrt{(b+1)^2+2\phi(x+1)}$ for all $\phi\in{\mathcal D}_1$. Then : 
\begin{equation*}
\Lambda_{t}(\phi,x)=\frac{1}{t}\log
E_{b(\phi,x)}(\exp(\{\phi+b-b(\phi,x))[F(Y_t)-F(y_0)-t/2]\})
\end{equation*}

When starting from $y_0 = 0$, the semi-group is deduced from that of $X$ :
\begin{align*}
\tilde{p}_{t}(0, y) = \sqrt{2\pi}K_{\alpha}\frac{e^{\gamma^2t/2}}{\sqrt{t}}\sum_{n \geq 0}\frac{\g(2n+ \alpha + \frac{3}{2})}{4^n n! \g(n+\alpha+1)}
(1-y^2)^{n+\alpha}  f_{T_1} \star f_{C_{2n+\gamma}}(1/t),
\end{align*}
where $p = -2(\alpha +1) = 2b \leq -2$  and $\gamma = -(p+1)/2 = \alpha +1/2$. Denote by
$${\mathcal D}=\{\phi \in \mathcal{D}_1 : G(\phi)=b + b(\phi,x)+\phi < 0\}\,.$$
For any $\phi\in{\mathcal D}$, the expectation
above is finite and a simple computation gives :
\begin{align*}
\Lambda_t(\phi,x) & =-\frac{\phi+b-b(\phi,x)}{2}+\frac{1}{t}\log\mathbb{E}_{b(\phi,x)}((1-Y_t^2)^{-(\phi+b-b(\phi,x))/2})
\\& = \Lambda(\phi,x) + \frac{1}{t}\log\frac{\sqrt{2\pi}K_{\alpha(\phi,x)}R_t(\phi,x)}{\sqrt t}
\end{align*} where
\begin{eqnarray*}
R_t(\phi,x) &=& \sum_{n \geq 0}\frac{\g(2n - b(\phi,x) + 1/2)}{4^nn! \g(n-b(\phi,x))}B\left(n - \frac{\phi+b +b(\phi,x)}{2},\frac{1}{2}\right) e^{\gamma^2t/2}f_{T_1} \star
f_{C_{2n+\gamma}}(\frac{1}{t})
\end{eqnarray*}and $B$ stands for the Beta function. Moreover, by dominated convergence theorem 
\begin{equation*}
\lim_{t \rightarrow \infty} \mathbb{E}_{b(\phi,x)}((1-Y_t^2)^{-(\phi+b-b(\phi,x))/2}) = \int_{-1}^{1} (1-y^2)^{-[\phi + b + b(\phi,x)]/2 \, -1} \, dy \, < \infty
\end{equation*}  
for $\phi \in \mathcal{D}$. Hence $\Lambda_t \rightarrow \Lambda$ as $t  \rightarrow \infty$. The following lemma details the domain ${\mathcal D}$ of $\Lambda_t$:
\begin{lem}\label{dom}
Denote by $$\phi_0=-\frac{(b+1)^2}{2(x+1)}\,.$$ \pn i)\underline{If $x<(b^2+3)/2(b-1)$}: then ${\mathcal
D}=(-\infty,\phi_0)$.
\medskip
\pn ii) \underline{If $(b^2+3)/2(b-1)<x<-1$}: then ${\mathcal D}=(-\infty,\phi_1)$ where $\phi_1$ is solution of $G(\phi)=0$.
\medskip
\pn iii) \underline{If $x>-1$}: then ${\mathcal D}=(\phi_0,\phi_1)$.
\end{lem}

In case i) of Lemma above, $\Lambda$ is steep. It achieves its unique minimum in $\phi_m$ solution of $$\partial
\Lambda/\partial \phi(\phi,x)=0\,,$$ i.e. $b(\phi,x)=x$. It is easy to see that
$$\phi_m=\frac{x+1}{2}-\frac{(b+1)^2}{2(x+1)} < \phi_0\,.$$
Hence, G\"artner-Ellis Theorem gives for $x < b < (b^2+3)/2(b-1)$,
\begin{eqnarray*}
\lim_{t\to\infty}\frac{1}{t}\log P(\hat b_t\leq x)= \lim_{t\to\infty}\frac{1}{t}\log P(S_{t,x} \leq 0)
 = \inf_{\phi \in ]\infty,\phi_0]}\Lambda(\phi,x) = \Lambda(\phi_m,x) = -\frac{1}{4}\frac{(x-b)^2}{x+1}.
\end{eqnarray*}
If $b<x< (b^2+3)/2(b-1)$, notice that $\phi_m>0$ and

\begin{eqnarray*}
\lim_{t\to\infty}\frac{1}{t}\log P(\hat b_t\geq x)= \lim_{t\to\infty}\frac{1}{t}\log P(S_{t,x} \geq 0)
 = \inf_{\phi \in (0,\phi_0]}\Lambda(\phi,x) = \Lambda(\phi_m,x) = -\frac{1}{4}\frac{(x-b)^2}{x+1}.
\end{eqnarray*}
\smallskip\par
 In cases ii) and iii) of Lemma \ref{dom}, $\Lambda$ is not steep. Nevertheless, if the infimum of $\Lambda$ is reached in $\overset{\circ}{\mathcal D}$,
 we can follow the scheme of Gartner--Ellis theorem for the change of probability in the infimum bound.
 This infimum is reached if and only if
\begin{equation}\label{dom1}\partial
\Lambda/\partial \phi(\phi_1,x) > 0.
\end{equation}
This above condition gives the following cases: denote by $x_0$ the unique solution $x<-1$ of $g(x) := 4x(x+1) -
(b-x)^2 = 0$. Since $g$ is decreasing on $]-\infty,-1]$ and $g(b^2+3/(2(b-1)) = (3/4)(b+1)^2 > 0 = g(x_0)$, then $x_0 > (b^2+3)/[2(b-1)]$.

\pn $\bullet$ if $(b^2+3)/2(b-1)<x<x_0<-1$, the derivative $\partial\Lambda/\partial\phi(\phi_1,x) > 0$, $\Lambda$ achieves its minimum on $\phi_m$ and
$$\lim_{t\to\infty}\frac{1}{t}\log P(\hat b_t \geq x)= \Lambda(\phi_m,x) = - \frac{(x-b)^2}{4(x+1)}\,.$$
\pn $\bullet$ if $x_0 < x < -1$ or $x > -1$, then $\partial\Lambda/\partial \phi(\phi_1,x) < 0$. We apply Theorem
\ref{pgd2} of the appendix, which is due to Zani \cite{Zani}. Let us verify that the assumptions are satisfied. Indeed,
the only singularity $\phi_1$ of $R_t$ comes from $B\left(n - [\phi+b +b(\phi,x)]/2,1/2\right)$ when $n=0$, and more
precisely, from $\Gamma(-[\phi+b + b(\phi,x)] / 2)\,.$ We can write
\begin{equation}
\Lambda_t(\phi,x)=\Lambda(\phi,x)+\frac{1}{t}\log\Gamma \left(-\frac{\phi+b+b(\phi,x)}{2}\right)+
\frac{1}{t}\log\frac{\sqrt{2\pi}K_{\alpha(\phi,x)}\tilde R_t(\phi,x)}{\sqrt t}\,,
\end{equation}
 where
\begin{equation}
\tilde{R_t}(\phi,x) = \frac{R_t(\phi,x)}{\Gamma(-[\phi+b+b(\phi,x)]/2)}
\end{equation}
Now $$\forall n \geq 0,\quad \frac{\d B\left(n -\frac{\phi+b + b(\phi,x)}{2},\frac{1}{2}\right)}{\d \Gamma \left(-\frac{\phi+b+b(\phi,x)}{2}\right)}$$ is analytic on some neighbourhood of $\phi_1$. Besides, $\phi_1$ is a pole of order one, i.e.
$$\lim_{\phi\to\phi_1\,,\phi<\phi_1}\frac{b+\phi+b(\phi,x)}{\phi-\phi_1}=c > 0\,,$$
and since $\lim_{\rho\to 0^+}\rho\Gamma(\rho)=1$, we can write
$$\frac{1}{t}\log\Gamma \left(-\frac{\phi+b+b(\phi,x)}{2}\right)=-\frac{\log(\phi_1-\phi)}{t}+\frac{h(\phi)}{t}\,.$$
The function $h$ is analytic on $\mathcal D$ and can be extended to an analytic function on $]\phi_1-\xi,\phi_1+\xi[$ for some positive $\xi$. Finally, we focus on $\tilde{R}_t(\phi,x)/\sqrt t$ and show that it  converges uniformly as $t \rightarrow \infty$. To proceed, we shall prove that this ratio is bounded from above and below away from $0$.
Setting  $A_n(t) : = e^{\gamma^2t/2}f_{T_1} \star f_{C_{2n+\gamma}}(1/t)$, one has :
\begin{align*}
\frac{A_n(t)}{\sqrt t} & \leq  \frac{e^{\gamma^2t/2}}{\sqrt{t}} \sum_{k,l \geq 0} U_{k,n} \int_0^{1/t}\exp-\frac{1}{2}\left[\frac{(2n+2k+\gamma)^2}{s}
+ \pi^2(l+\frac{1}{2})^2(\frac{1}{t} - s)\right] \frac{ds}{s^{3/2}}
\\& = \frac{e^{\gamma^2t/2}}{\sqrt{t}} \sum_{k,l \geq 0} U_{k,n} \int_t^{\infty}\exp-\frac{1}{2}\left[(2n+2k+\gamma)^2s
+ \pi^2(l+\frac{1}{2})^2(\frac{s-t}{ts})\right] \frac{ds}{\sqrt s} 
\\& < e^{-2n^2}\sum_{k,l \geq 0}U_{k,n}e^{-2k^2}\int_t^{\infty}\exp-\frac{1}{2}\left[(2n+2k+\gamma)^2(s-t) + \pi^2(l+\frac{1}{2})^2(\frac{s-t}{ts})\right] \frac{ds}{\sqrt{ts}}
\\&  = e^{-2n^2}\sum_{k,l \geq 0}U_{k,n}e^{-2k^2}\int_0^{\infty}\exp-\frac{1}{2}\left[(2n+2k+\gamma)^2 s + \pi^2 l^2(\frac{s}{t(t+s)})\right] \frac{ds}{\sqrt{t(t+s)}}
\end{align*}
with
$$U_{k,n} = \frac{\Gamma(2n+k + \gamma)2^{2n+\gamma}(2n+2k+\gamma)}{k!\Gamma(2n+\gamma)}\,.$$
Let $\Theta(x) = \sum_{l \in \mathbb{Z}}e^{-\pi l^2 x} = 1 +2\sum_{l \geq 1}e^{-\pi l^2 x}$ denote the Jacobi Theta function. Then

\begin{equation*}
\frac{A_n(t)}{\sqrt{t}} <  e^{-2n^2}\left[\sum_{k \geq 0}U_{k,n}e^{-2k^2}\int_0^{\infty}
\exp-\left[\frac{(2n+2k+\gamma)^2 s}{2}\right] \Theta\left(\frac{\pi s}{2t(t+s)}\right)\frac{ds}{\sqrt{t(t+s)}} + C(n,t)\right]
\end{equation*}
where 
\begin{equation*}
C(n,t) =  \frac{1}{2\sqrt t}\sum_{k,l \geq 0}U_{k,n}e^{-2k^2}\int_0^{\infty}\exp-\left[\frac{(2n+2k+\gamma)^2 s}{2}\right] \frac{ds}{\sqrt{t+s}}
\end{equation*}

Recall that $\Theta(x) = (1/\sqrt{x})\Theta(1/x)$, which yields :

\begin{equation*}
\frac{A_n(t)}{\sqrt t} < e^{-2n^2}\sum_{k \geq 0}U_{k,n}e^{-2k^2}\int_0^{\infty}\exp-\left[\frac{(2n+2k+\gamma)^2s}{2}\right] \Theta\left(\frac{2t(t+s)}{\pi s}\right)\frac{ds}{\sqrt s}
+ \frac{C(n)}{2\sqrt t}
\end{equation*}
where 
\begin{equation*}
C(n) = e^{-2n^2}\sum_{k,l \geq 0}U_{k,n}e^{-2k^2}\int_0^{\infty}\exp-\left[\frac{(2n+2k+\gamma)^2 s}{2}\right] \frac{ds}{\sqrt s}
\end{equation*}

Since $e^{-l^2z} < e^{-lz}$, then $\Theta(z)\leq 3$ for $z > 1$. Hence, as $2t/\pi \leq 2t(t+s)/(\pi s)$, then for $t$ large enough : 
\begin{align*}
\frac{A_n(t)}{\sqrt t} & < 3 e^{-2n^2}\sum_{k \geq 0}U_{k,n}e^{-2k^2}\int_0^{\infty}\exp-\left[\frac{(2n+2k+\gamma)^2 s}{2}\right] \frac{ds}{\sqrt s} + C(n) 
\end{align*}
This gives a lower bound for $\tilde{R}_t/\sqrt t$. Besides,
\begin{align*}
&\frac{\tilde{R}_t(\phi,x)}{\sqrt t}  > \frac{\sqrt{\pi}\g(1/2-b(\phi,x))}{\Gamma(-b(\phi,x))\Gamma\{[1- (\phi+b+b(\phi,x)]/2\}} \frac{A_0(t)}{\sqrt t}
\\& = C(b,\phi,x) \sum_{k,l \geq 0} (-1)^kV_k\int_0^{\infty}\exp-\frac{1}{2}\left[(2k+\gamma)^2 s + \pi^2(l+\frac{1}{2})^2(\frac{s}{t(t+s)})\right] \frac{ds}{\sqrt{t(t+s)}}
\end{align*}
where $V_k(t) := U_{k,0}e^{-2k(k+\gamma)t}$. One may choose $t$ large enough independent of $k$ such that $V_k(t) \geq V_{k+1}(t)$ for all $k \geq 0$. 
In fact, such $t$ satisfies :
\begin{equation*}
e^{2(2k+\gamma+1)t} \geq e^{2t} \geq \sup_{k \geq 0}\frac{U_{k+1,0}}{U_{k,0}} = \sup_{k \geq 0}\frac{(k+\gamma)(2k + \gamma +2)}{(k+1)(2k+\gamma)}
\end{equation*}
Then :
\begin{align*}
\frac{\tilde{R}_t}{\sqrt t} &> C(b,\phi,x) [V_0(t) - V_1(t)] \sum_{l \geq 0} \int_0^{\infty}\exp-\frac{1}{2}\left[\gamma^2 s + \pi^2(l+\frac{1}{2})^2(\frac{s}{t(t+s)})\right] \frac{ds}{\sqrt{t(t+s)}}
 \\&> C(b,\phi,x) [\gamma 2^{\gamma} - V_1(t)] \sum_{l \geq 0} \int_0^{\infty}\exp-\frac{1}{2}\left[\gamma^2 s + \pi^2(l+1)^2(\frac{s}{t(t+s)})\right] \frac{ds}{\sqrt{t(t+s)}}
\\& = \frac{C(b,\phi,x)}{2} [\gamma 2^{\gamma} - V_1(t)] \left\{\int_0^{\infty}e^{-\gamma^2 s/2}\Theta\left(\frac{\pi s}{2t(t+s)}\right) \frac{ds}{\sqrt{t(t+s)}} - C(t)\right\}.
\end{align*}
where 
\begin{equation*}
C(t) = \frac{1}{2\sqrt t}\int_0^{\infty}e^{-\gamma^2 s/2} \frac{ds}{\sqrt{(t+s)}} < c \int_0^{\infty}e^{-\gamma^2 s/2} \frac{ds}{\sqrt{s}}, \quad c < \sqrt{\frac{2}{\pi}}  .
\end{equation*}
for $t$ large enough. Following the same scheme as for the upper bound, one gets : 
\begin{align*}
\frac{\tilde{R}_t}{\sqrt t} & > \frac{C(b,\phi,x)}{2} \gamma 2^{\gamma} \left\{\sqrt{\frac{2}{\pi}}\int_0^{\infty}e^{-\gamma^2 s/2}\Theta\left(\frac{2t(t+s)}{\pi s}\right) \frac{ds}{\sqrt s}
 - C(t)\right\}
\\& >  \frac{C(b,\phi,x)}{2} \gamma 2^{\gamma} \left(\sqrt{\frac{2}{\pi}}- c\right) \int_0^{\infty}e^{-\gamma^2 s/2} \frac{ds}{\sqrt s} > 0. \qquad \qquad \square
\end{align*}
As a result,
$$\lim_{t\to\infty}\frac{1}{t}\log P(\hat b_t \geq x)=\Lambda(\phi_1,x)=-(x+2+\sqrt{(b-x)^2+4(x+1)}),$$
which ends the proof of Theorem \ref{pgd1}. $\hfill \blacksquare$

\subsection{Jacobi-squared Bessel processes duality.}
By It\^o's formula and L\'evy criterion, one claims that $(Y_t^2)_{t \geq 0}$ is a Jacobi process of parameters $d=1,\,d' = -2b \geq 2$. Indeed :
\begin{align*}
dZ_t : = d(Y_t^2) & = 2Y_tdY_t + \langle Y \rangle_t = 2Y_t\sqrt{1-Y_t^2}dW_t + [(2b-1)Y_t^2 + 1] dt
\\& = 2\sqrt{Z_t(1-Z_t)}\textrm{sgn}(Y_t)dW_t  + [(2b-1)Z_t + 1] dt
\\& = 2\sqrt{Z_t(1-Z_t)}dB_t  + [(2b-1)Z_t + 1] dt
\end{align*}
Using the skew product previously stated, there exists $R$, a squared Bessel process of dimension $d' = 2(\nu+1) = -2b$ and starting from $r$ so that :
\begin{equation*}
\hat{\nu}_t : = - \hat{b}_t - 1 = \frac{\log(1-Y_t^2) + t}{2\int_0^t \frac{Y_s^2}{1-Y_s^2}ds}
\end{equation*}
is another estimator of $\nu$ based on a Jacobi trajectory observed till time $t$. Set $t = \log u$, then
\begin{equation*}
 \hat{\nu}_{\log u}^1 = \frac{\log[u(1-Y_{\log u}^2)]}{2\int_0^{\log u} \frac{Y_s^2}{1-Y_s^2}ds} = \frac{\log[u(1-Y_{\log u}^2)]}{2\int_1^{u} \frac{Y_{\log s}^2}{s(1-Y_{\log s}^2)}ds}
\end{equation*}
and $ \{\hat{\nu}_{\log u}^1\}_u$ satisfies a LDP with speed $\log u$ and rate function $J_{-(\nu+1)}(-(x+1))$.

When starting at $R_0 = 1$, the MLE of $\nu$ based on a Bessel trajectory is given by (cf \cite{Zani}, p. 132) :
\begin{equation*}
\hat{\nu}_t^1 = \frac{\int_0^t\frac{dX_s}{X_s} - 2\int_0^t \frac{ds}{X_s}}{2\int_0^t \frac{ds}{X_s}} = \frac{\log(X_t)}{2\int_0^t\frac{ds}{X_s}}
\end{equation*}
with associated rate function :
\begin{equation*}
I_{\nu}(x) =  \left\{\begin{array}{lccr}
\frac{(x-\nu)^2}{4x} & \textrm{if} & x \geq x_1    & : = \frac{-(\nu+2) + 2\sqrt{\nu^2+\nu+1}}{3} \\
1-x + \sqrt{(\nu-x)^2 -4x} & \textrm{if} & x < x_1& \end{array}
\right.
\end{equation*}
 A glance at both rate functions gives $I_{\nu}(x) = J_{-(\nu+1)}(-(x+1))$ and $x_0 = -(x_1+1)$.

\section{Appendix}
Let $\{Y_t\}_{t\geq 0}$ be a family of real random variables defined on $(\Omega,{\mathcal F},P)$, and denote by
$\mu_t$ the distribution of $Y_t$. Suppose $-\infty < m_t := EY_t < 0$. We look for large deviations bounds for
$P(Y_t\geq y)$. Let $\Lambda_t$ be the n.c.g.f. of $Y_t$:
$$\Lambda_t(\phi)=\frac{1}{t}\log\,E(\exp\{\phi tY_t\})\,,$$
and denote by $D_t$ the domain of $\Lambda_t$. We assume that there exists $0< \phi_1 < \infty$ such that for any $t$
$$\sup \{ \phi : \phi \in D_t \} = \phi_1$$
and $[0 , \phi_1 ) \subset D_t$. We assume also that for $\phi \in D$
\begin{assum}
\label{hyp3}
\begin{equation}\label{ache}
\Lambda_{t} (\phi) = \Lambda(\phi) - \frac{\alpha}{t} \log (\phi_{1} - \phi) +\frac{R_{t} (\phi)}{t}
\end{equation}
where \pn $\bullet$ $\alpha>0$ \pn $\bullet$ $\Lambda$ is analytic on $(0,\phi_1)$, convex, with finite limits at
endpoints, such that $\Lambda'(0) < 0$, $\Lambda'(\phi_{1}) < \infty$, and $\Lambda''(\phi_1)>0$. \pn $\bullet$ $R_{t}$
is analytic on $(0, \phi_1 )$ and admits an analytic extension on a strip $D_{\beta}^{\gamma}=(\phi_1-\beta,\phi_1 +
\beta) \times(-\gamma,\gamma)$, where $\beta$ and $\gamma$ are independent of $t$. \pn $\bullet$ $R_{t} (\phi)$
converges as $t \rightarrow \infty $ to some $R(\phi )$ uniformly on any compact of $D_{\beta}^{\gamma}$.
\end{assum}

\begin{teo}
\label{pgd2}
Under \ref{hyp3} \pn For any $\Lambda'(0) < y<\Lambda'(\phi_1)$,
\begin{equation}
\label{rate} \lim_{t\rightarrow +\infty}\frac{1}{t}\log\,P(Y_t\geq y)=
-\sup_{\phi\in(0,\phi_1)}\{y\phi-\Lambda(\phi)\}\,.
\end{equation}
For any $y\geq\Lambda'(\phi_1)$,
\begin{equation}
\label{ratelin} \lim_{t\rightarrow +\infty}\frac{1}{t}\log\,P(Y_t\geq y)= -y\phi_1+\Lambda(\phi_1)\,.
\end{equation}
The rate function is continuously differentiable with a linear part.
\end{teo}

\end{document}